\documentclass{article}
\usepackage[french,british]{babel}

\usepackage[T1]{fontenc}

\usepackage{amsmath}
\usepackage{amsfonts}
\usepackage{amsbsy}
\usepackage{amssymb}

\usepackage{amsthm}

\theoremstyle{plain}

\newtheorem{lemma}{Lemma}
\newtheorem{theorem}{Theorem}

\newtheorem{corollary}{Corollary}

\newtheorem{definition}{Definition}

\newtheorem{Proposition}{Proposition}[section]

\newtheorem{exm}{Example}

\begin{document}

\selectlanguage{british}

\title{Minimax extrapolation problem for periodically correlated stochastic sequences with missing observations}

\author{
Iryna Golichenko\thanks{Department of Mathematical Analysis and Probability Theory, National Technical  University of Ukraine ''Igor Sikorsky Kyiv Politechnic Institute'',
 Kyiv 03056, Ukraine,\newline e-mail: idubovetska@gmail.com},
Oleksandr Masyutka\thanks{Department of Mathematics and Theoretical Radiophysics, Taras Shevchenko National University of Kyiv, 01601, Ukraine, e-mail: omasyutka@gmail.com},
Mikhail Moklyachuk\thanks
{Department of Probability Theory, Statistics and Actuarial
Mathematics, Taras Shevchenko National University of Kyiv, Kyiv 01601, Ukraine, e-mail: moklyachuk@gmail.com}
}

\date{\today}

\date{\today}

\maketitle

\renewcommand{\abstractname}{Abstract}
\begin{abstract}
 The problem of optimal estimation of the linear functionals
 which depend on the
unknown values of a periodically correlated stochastic sequence
${\zeta}(j)$ from observations of the sequence
${\zeta}(j)+{\theta}(j)$ at points   $j\in\{\dots,-n,\dots,-2,-1,0\}\setminus S$, $S=\bigcup _{l=1}^{s-1}\{-M_l\cdot T+1,\dots,-M_{l-1}\cdot T-N_{l}\cdot T\}$, is
considered, where ${\theta}(j)$ is an uncorrelated
with ${\zeta}(j)$ periodically correlated stochastic sequence.
Formulas for calculation the mean square error and the spectral
characteristic of the optimal estimate of the functional
$A\zeta$ are proposed in the case where spectral densities of the sequences are
exactly known.
Formulas that determine the least favorable spectral densities and the minimax-robust spectral characteristics of the optimal estimates of functionals are proposed in the case of spectral uncertainty, where the spectral densities are not exactly known while some sets of admissible spectral densities are specified.
\end{abstract}

\vspace{2ex}

\maketitle

\vspace{2ex}

\section{Introduction}

The problem of estimation of the unknown values of stochastic sequences and processes is of constant interest in the theory of stochastic processes. The formulation of the interpolation, extrapolation and filtering problems for stationary stochastic sequences with known spectral densities and reducing them to the corresponding problems of the theory of functions belongs to
Kolmogorov (see, for example, selected works by Kolmogorov, 1992). Effective methods of solution of the estimation problems for stationary stochastic sequences and processes were developed by Wiener (1966) and Yaglom (1987). Further results are presented in the books  by Rozanov (1967) and Hannan (1970).

In 1958  Bennett Introduced the notion of cyclostationarity as a phenomenon and property of the process, which describes signals in channels of communication. Studying the statistical characteristics of information transmission, he calls the group of telegraph signals the cyclostationary process, that is the process whose group of statistics changes periodically with time.
Gardner \& Franks (1975)  highlights the greatest similarity of  cyclostationary processes, which are a subclass of nonstationary processes, with stationary processes.
Gardner (1994) presented the bibliography of works   in which properties and applications of cyclostationary processes were investigated.
Recent developments and applications of cyclostationary signal analysis are reviewed in the papers by Gardner, Napolitano \& Paura (2006) and Napolitano  (2016).
Note, that in other sources cyclostationary processes are called periodically stationary, periodically nonstationary, periodically correlated. We will use the term periodically correlated processes.
 Gladyshev  (1961) was the first who started the analysis of spectral properties and representation of periodically correlated sequences based on its connection with vector stationary sequences. He formulated the necessary and sufficient conditions for determining of periodically correlated sequence in terms of the correlation function.
Makagon (1999), Makagon et al. (2011) presented  detailed spectral analysis of periodically correlated sequences.
Main ideas of the research  of periodically correlated sequences are outlined in the book by Hurd \& Miamee (2007).

Since stochastic processes often accompanied with undesirable noise it is naturally to assume that the exact value of  spectral density is unknown and the model of process  is given by a set of restrictions on spectral density.
Vastola \& Poor  (1983) have demonstrated that the described procedure
can result in significant increasing of the value of error. This is
a reason for searching estimates which are optimal for all densities
from a certain class of admissible spectral densities. These
estimates are called minimax since they minimize the maximal value
of the error of estimates. A survey of results in minimax (robust)
methods of data processing can be found in the paper by Kassam \& Poor (1985).

Grenander (1957) was the first who
proposed the minimax approach to the extrapolation problem for
stationary processes. Formulation and investigation of the problems of extrapolation, interpolation and filtering of linear functionals which depend on the unknown values of stationary sequences and processes  from observations with and without noise are presented by
Moklyachuk (2008), (2015).
Results of investigation of the problems of optimal estimation of vector-valued stationary sequences and processes are published by
Moklyachuk \& Masyutka  (2006 -- 2012).
In their book
Luz \& Moklyachuk (2019) presented results of investigation of the minimax estimation problems for linear functionals which depends on unknown
values of stochastic sequence with stationary increments.
 Golichenko \& Moklyachuk  (2012 -- 2016) investigated the
interpolation, extrapolation and filtering problems of linear functionals from periodically correlated stochastic sequences and processes.
  The interpolation and filtering problems
for stationary sequences with missing values was examined by Moklyachuk, Masyutka \& Sidei (2019).
The interpolation problem of linear functionals from periodically correlated stochastic sequences with missing observations was investigated by
Golichenko \& Moklyachuk  in (2020).

In this paper we presented results of investigation of the problem of optimal linear estimation of
the functional
$A{\zeta}=\sum_{j=1}^{\infty}{a}(j){\zeta}(j),$
 which depends on the
unknown values of a periodically correlated stochastic sequence
${\zeta}(j)$ from observations of the sequence
${\zeta}(j)+{\theta}(j)$ at points $j\in\{\dots,-n,\dots,-2,-1,0\}\setminus S$, $S=\bigcup _{l=1}^{s-1}\{-M_l\cdot T+1,\dots,-M_{l-1}\cdot T-N_{l}\cdot T\}$, where
${\theta}(j)$ is an uncorrelated with ${\zeta}(j)$ periodically correlated stochastic sequence.
Formulas for calculation of the mean square error and the spectral
characteristic of the optimal estimate of the functional
$A\zeta$ are proposed in the case where spectral densities are
exactly known. Formulas that determine the least favorable spectral densities and the minimax-robust spectral characteristics of the optimal estimates of functionals are proposed in the case of spectral uncertainty, where the spectral densities are not exactly known while some sets of admissible spectral densities are specified.

\section{Periodically correlated and multidimensional stationary  sequences}

The  term \textit {periodically correlated} process was introduced by Gladyshev (1961) while Bennett (1958) called random and periodic processes
\textit {cyclostationary} process.

Periodically correlated sequences are stochastic sequences that have
periodic structure (see the book by Hurd \& Miamee (2007).

\begin{definition} A complex valued stochastic
sequence ${\zeta}(n),n\in{\mathbb Z}$ with zero mean,
$\textsf{E}{\zeta}(n)=0$, and finite variance, $\textsf{E}|{\zeta}(n)|^{2}<+\infty$,
is called cyclostationary or periodically
correlated (PC) with period $T$ ($T$-PC) if for every
$n,m\in{\mathbb Z}$
\begin{equation}
\label{2.1}
\textsf{E}{\zeta}(n+T)\overline{\zeta(m+T)}=R(n+T,m+T)=R(n,m)
\end{equation}
and there are no smaller values of $T>0$ for
which (1)  holds true.
\end{definition}

\begin{definition} {A complex valued
T-variate stochastic  sequence $\vec{\xi}(n)=\left\{ \xi_{\nu}(n)
\right\}_{\nu = 1 }^{T},$ $n\in {\mathbb Z}$ with zero mean,
$\textsf{E}\xi_{\nu}(n)=0, \nu=1,\dots,T$, and $\textsf{E}||\vec {\xi}(n)||^{2}<\infty$
is called stationary if for all $n ,m\in {\mathbb Z}$ and
$\nu,\mu\in \{1,\dots,T\}$}
\[
\textsf{E}\xi_{\nu}(n)\overline{\xi_{\mu}(m)}=R_{\nu\mu}(n,m)=R_{\nu\mu}(n-m).
\]
\noindent If this is the case, we denote $R(n)=\left\{ R_{\nu\mu}(n)
\right\}_{\nu,\mu = 1 }^{T}$ and call it the \textit {covariance
matrix} of T-variate stochastic sequence $\vec{\xi}(n)$.
\end{definition}

\begin{Proposition} {\rm{(Gladyshev, 1961)}}.
A stochastic sequence $\zeta(n)$ is PC with period $T$ if and only if there exists
a $T$-variate stationary sequence $\vec {\xi}(n)=\left\{ \xi_{\nu}(n)
\right\}_{\nu = 1 }^{T}$ such that $\zeta(n)$ has the representation
\begin{equation}
\label{2.2}
\zeta(n)=\sum_{\nu=1}^{T}{e}^{2\pi in\nu/T}\xi_{\nu}(n),\,\,n\in {\mathbb
Z}.
\end{equation}
\noindent The sequence $\vec \xi(n)$ is called  \textit {generating sequence} of the sequence $\zeta(n)$.
\end{Proposition}

\begin{Proposition}{\label{1}}{\rm{(Gladyshev, 1961)}}. A complex valued stochastic
sequence ${\zeta}(n),n\in{\mathbb Z}$ with zero mean and finite variance is PC with period $T$ if and only if the  $T$-variate blocked sequence $\vec \zeta(n)$ of the form
\begin{equation}
\label{block}
[\vec{\zeta}(n)]_{p}=\zeta(nT+p),\,\,n\in {\mathbb
Z},p=1,\dots,T
\end{equation}
 is stationary.
 \end{Proposition}

We will denote by $f^{\vec \zeta}(\lambda)=\left\{ f^{\vec \zeta}_{\nu\mu}(\lambda)
\right\}_{\nu,\mu = 1 }^{T}$ the matrix valued
spectral density function of the $T$-variate  stationary
sequence $\vec \zeta(n)=(\zeta_1(n),\dots,\zeta_T(n))^{\top}$ arising from the $T$-blocking (3) of a univariate T-PC
sequence $\zeta(n)$.

\section{The classical projection method of linear extrapolation}

Let $\zeta(j)$ and $\theta(j)$ be uncorrelated T-PC stochastic
sequences. Consider the problem of optimal linear estimation  of the
functional $$A{\zeta}=\sum_{j=1}^{\infty}{a}(j){\zeta}(j),$$
 that depends on the unknown values of T-PC stochastic
sequence $\zeta(j)$, based on
observations of the sequence $\zeta(j)+\theta(j)$ at points $j\in \{...,-n,...,-1,0\} \setminus S,$ $S=\bigcup _{l=1}^{s}\{-M_l\cdot T+1,\dots,-M_{l-1}\cdot T-N_{l}\cdot T\},\,
 M_l=\sum_{k=0}^l(N_k+K_k),\,\,N_0=K_0=0$.

Let assume that the  coefficients $a(j), j\geq1$  which determine the functional $A\zeta$ satisfy condition
\begin{equation}
\label{suma}
\sum_{j=1}^{\infty} |a(j)|<\infty
\end{equation}
and are of the form
\begin{equation}
\label{aj}
a(j)=a\left(\left(j-\left[\frac{j}{T}\right]T\right)+\left[\frac{j}{T}\right]T\right)=a(\nu+\tilde{j}T)=a(\tilde{j})e^{2\pi i\tilde{j}\nu/T},\,
\end{equation}
$$
\nu=1,\dots,T,\, \tilde{j}\geq0,$$
where  $\nu=T$ and $\tilde{j}=\lambda-1$, if $j=T\cdot \lambda, \,\lambda \in \mathbb{Z},$ or
$$
 a(j)=a(T\cdot \lambda)=a(T+(\lambda-1)T)=a(\lambda-1)e^{2\pi i(\lambda-1)T/T}.
$$
Under the condition (4)  the functional $A\zeta$ has the finite second moment.

Using Proposition 2.2,  the
 linear functional $A{\zeta}$ can be written as follows
$$
A\zeta=\sum_{j=1}^{\infty}{a}(j){\zeta}(j)=\sum_{{\widetilde{j}}=0}^{\infty}{a}(\widetilde{j})\sum_{\nu=1}^{T}e^{2\pi
i\widetilde{j}\nu/T}\zeta(\nu+\widetilde{j}T)=
$$
$$
=\sum_{\widetilde{j}=0}^{\infty}\sum_{\nu=1}^{T}{a}(\widetilde{j})e^{2\pi
i\widetilde{j}\nu/T}\zeta_{\nu}(j)=\sum_{\widetilde{j}=0}^{\infty}\vec {a}^{\top}(\widetilde{j})\vec
{\zeta}(\widetilde{j})=A\vec {\zeta},
$$
\noindent where
\begin{equation}
\label{avec}
\vec {a}^{\top}(\widetilde{j})=\left(
a_{1}(\widetilde{j}),
\dots,
a_{T}(\widetilde{j})
\right),\,a_{\nu}(\widetilde{j})={a}(\widetilde{j})e^{2\pi
i\widetilde{j}\nu/T},\,\nu=1,\dots,T,\end{equation}
$\vec {\zeta}(\widetilde{j})=\left\{ \zeta_{\nu}(\widetilde{j})
\right\}_{\nu = 1}^{T}$ is $T$-variate stationary sequence, obtained by the $T$-blocking (3) of univariate $T$-PC sequence $\zeta(j), \, j \geq1$.

Let $\vec {\zeta}(j)$ and $\vec {\theta}(j)$ be uncorrelated T-variate
stationary stochastic sequences with the spectral density matrices
$f^{\vec \zeta}(\lambda)=\left\{ f^{\vec \zeta}_{\nu\mu}(\lambda)
\right\}_{\nu,\mu = 1 }^{T}$ and $f^{\vec \theta}(\lambda)=\left\{
f^{\vec \theta}_{\nu\mu}(\lambda) \right\}_{\nu,\mu= 1}^{T}$,
respectively. Consider the problem of optimal linear estimation of
the functional $$A\vec {\zeta}=\sum_{\widetilde{j}=0}^{\infty}\vec {a}^{\top}(\widetilde{j})\vec
{\zeta}(\widetilde{j}),$$ that depends on the unknown values of sequence  $\vec {\zeta}(\widetilde{j})$,
based on observations of the sequence $\vec
{\zeta}(\widetilde{j})+\vec{\theta}(\widetilde{j})$ at points $\widetilde{j}\in\{...,-n,...,-1\}
\setminus \widetilde{S},$ $\widetilde{S}=\bigcup _{l=1}^{s}\{-M_l,\dots,-M_{l-1}-N_{l}-1\},\,
 M_l=\sum_{k=0}^l(N_k+K_k),\,\,N_0=K_0=0,$.

Let the spectral densities $f^{\vec \zeta}(\lambda)$ and $f^{\vec
\theta}(\lambda)$ satisfy the minimality condition
\begin{equation}
\label{3.5}
\int_{-\pi}^{\pi}{Tr{\left[ {(f^{\vec \zeta}(\lambda)+f^{\vec
\theta}(\lambda))^{-1}} \right]}} d\lambda <+{\infty}.
\end{equation}
\noindent Condition (7) is necessary and sufficient in order that
the error-free extrapolation of unknown values of the sequence $\vec
{\zeta}(j)+\vec {\theta}(j)$ is impossible \cite{Rozanov}.

Denote by $L_{2}(f)$
the Hilbert space of vector valued functions $\vec  b(\lambda)=\left\{
b_{\nu}(\lambda) \right\}_{\nu = 1 }^{T}$ that are integrable with
respect to a measure with the density $f(\lambda)=\left\{
f_{\nu\mu}(\lambda) \right\}_{\nu,\mu = 1 }^{T}$:
$$
\int_{-\pi}^{\pi}\vec b^{\top}(\lambda)f(\lambda)\overline{\vec b(\lambda)}d\lambda=\int_{-\pi}^{\pi}\sum_{\nu,\mu=1}^{T}b_{\nu}(\lambda) f_{\nu\mu}(\lambda) \overline{b_{\mu}(\lambda)}d\lambda<+\infty.
$$

 Denote by $L_{2}^{s}(f)$ the subspace in $L_{2}(f)$
generated by functions $$e^{i\widetilde{j}\lambda}\delta_{\nu},\delta_{\nu}=\left\{
\delta_{\nu\mu} \right\}_{\mu =1}^{T}, \,\nu=1,\dots,T,\,\widetilde{j}\in \{...,-n,...,-1\}\setminus \widetilde{S},$$ where $\delta_{\nu\nu}=1,\delta_{\nu\mu}=0$
for $\nu\neq\mu$.

Every linear estimate $\widehat{A\vec {\zeta}}$ of the functional
$A\vec {\zeta}$ from observations of the sequence $\vec
{\zeta}(\widetilde{j})+\vec{\theta}(\widetilde{j})$ at points $\widetilde{j}\in \{...,-n,...,-1\}
\setminus \widetilde{S}$ has the form
\begin{equation}
\label{est}
\widehat{A\vec \zeta}=\int_{-\pi}^{\pi}\vec h^{\top}(e^{i\lambda})(Z^{\vec \zeta}(d\lambda)+Z^{\vec \theta}(d\lambda))=\int_{-\pi}^{\pi}\sum_{\nu=1}^{T}h_{\nu}(e^{i\lambda})(Z_{\nu}^{\vec \zeta}(d\lambda)+Z_{\nu}^{\vec \theta}(d\lambda)),
\end{equation}
\noindent where $Z^{\vec \zeta}(\Delta)=\left\{ Z_{\nu}^{\vec \xi}(\Delta)
\right\}_{\nu =1 }^{T}$ and $Z^{\vec \theta}(\Delta)=\left\{
Z_{ \nu}^{\vec \eta}(\Delta) \right\}_{\nu = 1 }^{T}$ are orthogonal random
measures of the sequences $\vec \zeta(\widetilde{j})$ and $\vec \theta(\widetilde{j})$, and
$\vec h(e^{i\lambda})=\left\{ h_{\nu}(e^{i\lambda}) \right\}_{\nu = 1}^{T}$ is the spectral characteristic of the estimate
$\widehat{A\vec \zeta}$. The function $\vec h(e^{i\lambda})\in
L_{2}^{s}(f^{\vec \zeta}+f^{\vec \theta})$.

The mean square error $\Delta(\vec h;f^{\vec \zeta},f^{\vec \theta})$ of the
estimate $\widehat{A\vec \zeta}$ is calculated by the formula
$$
\Delta(\vec h;f^{\vec \zeta},f^{\vec \theta})=E|A\vec
{\zeta}-\widehat{A\vec \zeta}|^{2}=
$$
\begin{equation}\label{3.55}
=\frac{1}{2\pi}\int_{-\pi}^{\pi}{\left[ {
A(e^{i\lambda})-\vec h(e^{i\lambda})} \right]}^{\top}f^{\vec
\zeta}(\lambda){\overline{\left[ {
A(e^{i\lambda})-\vec h(e^{i\lambda})} \right]}}
d\lambda+\end{equation}
$$+\frac{1}{2\pi}\int_{-\pi}^{\pi} \vec h^{\top}(e^{i\lambda})f^{\vec
\theta}(\lambda)\overline{\vec h(e^{i\lambda})}d\lambda,
$$
$$
A(e^{i\lambda})=\sum_{\widetilde{j}=0}^{\infty}\vec {a}(\widetilde{j})e^{i\widetilde{j}\lambda}.
$$
\noindent The spectral characteristic $\vec h(f^{\vec \zeta},f^{\vec
\theta})$ of the optimal linear estimate of $A\vec {\zeta}$
minimizes the mean square error
\begin{equation}
\label{3.6}
\Delta(f^{\vec \zeta},f^{\vec \theta})=\Delta(\vec h(f^{\vec \zeta},f^{\vec
\theta});f^{\vec \zeta},f^{\vec \theta})=\mathop {\min }\limits_{\vec h \in
L_{2}^{s}(f^{\vec \zeta}+f^{\vec \theta})} \Delta (\vec h;f^{\vec
\zeta},f^{\vec \theta})=\mathop {\min }\limits_{\widehat{A\vec
\zeta}}E|A\vec {\zeta}-\widehat{A\vec \zeta}|^{2}.
\end{equation}
\noindent With the help of the
Hilbert space projection method proposed by Kolmogorov we can
find a
solution of the optimization problem (10). The optimal linear estimate $\widehat{A\vec \zeta}$ is a projection of the functional $A\vec \zeta$ on the subspace $H^{s}[\vec \zeta+\vec \theta]=H^{s}[\zeta_\nu(\widetilde{j})+\theta_\nu(\widetilde{j}),\, \widetilde{j}\in \{...,-n,...,-1\}\backslash \widetilde{S}, \nu=1,\dots,T]$ of the Hilbert space $H=\{\zeta: \textsf{E}\zeta=0,\, \textsf{E}|\zeta|^2<\infty\}$,
 generated by  values $\zeta_\nu(\widetilde{j})+\theta_\nu(\widetilde{j}),\, \widetilde{j}\in\{...,-n,...,-1\}\backslash \widetilde{S}, \nu=1,\dots,T$. The projection is characterized by following conditions

1) $\widehat {A \vec \zeta} \in H^{s}[\vec \zeta+\vec \theta],$

 2) $A \vec \zeta-\widehat {A\vec \zeta}\perp H^s[\vec \zeta+\vec \theta].$

The condition 2) gives us the possibility to derive the  formula for  spectral characteristic of the estimate
$$
\vec h^{\top}(f^{\vec \zeta},f^{\vec \theta})= \left( A^{\top}
(e^{i\lambda} )f^{\vec \zeta}(\lambda) - C^{\top} (e^{i\lambda
})\right)\, \left[ f^{\vec \zeta}(\lambda)+f^{\vec \theta}(\lambda)
\right]^{-1}=
$$
\begin{equation}
\label{3.7}
=A^{\top}(e^{i\lambda})-\left(
{A^{\top}(e^{i\lambda})f^{\vec
\theta}(\lambda)+C^{\top}(e^{i\lambda})}\right) \left[ f^{\vec
\zeta}(\lambda) +f^{\vec \theta}(\lambda)\right]^{-1},
\end{equation}
\noindent where
$$
C(e^{i\lambda})=\sum_{n\in \Gamma} \vec
{c}(n)e^{in\lambda},$$
where $\Gamma=\widetilde{S}\cup\{0,1,2,...\}$ and
$\vec
{c}(n), n\in \Gamma,$ are unknown vectors of coefficients.

Condition 1) is satisfied if the system of equalities
\begin{equation}
\label{syst}
\int_{-\pi}^{\pi}\vec h(f^{\vec \zeta},f^{\vec
\theta})e^{-im\lambda}d\lambda=0, m\in \Gamma
\end{equation}
holds true.

The last equalities (12) provide the following relations
$$\sum_{\widetilde{j}=0}^\infty \vec a^{\top}(\widetilde{j})\frac{1}{2\pi}\int_{-\pi}^{\pi}f^{\vec \zeta}(\lambda)(f^{\vec \zeta}(\lambda)+f^{\vec \theta}(\lambda))^{-1}e^{i\lambda(\widetilde{j}-m)}d\lambda=
$$
\begin{equation}
\label{2}
\sum_{n\in \Gamma} \vec c^{\top}(n)\frac{1}{2\pi}\int_{-\pi}^{\pi}(f^{\vec \zeta}(\lambda)+f^{\vec \theta}(\lambda))^{-1}e^{i\lambda(n-m)}d\lambda,\, \forall m\in \Gamma.\end{equation}

Denote the Fourier coefficients of the matrix functions $(f^{\vec \zeta}(\lambda)+f^{\vec \theta}(\lambda))^{-1}$  and $f^{\vec \zeta}(\lambda)(f^{\vec \zeta}(\lambda)+f^{\vec \theta}(\lambda))^{-1}$ as
$$
B(m-n)=\frac{1}{2\pi}\int_{-\pi}^{\pi}(f^{\vec \zeta}(\lambda)+f^{\vec \theta}(\lambda))^{-1}e^{i\lambda(n-m)}d\lambda,
$$
$$
R(m-\widetilde{j})=\frac{1}{2\pi}\int_{-\pi}^{\pi}f^{\vec \zeta}(\lambda)(f^{\vec \zeta}(\lambda)+f^{\vec \theta}(\lambda))^{-1}e^{i\lambda(\widetilde{j}-m)}d\lambda,
$$
$$
n,m \in \Gamma,\, \widetilde{j}=0,1,2....
$$

Denote by $\vec {\mathbf a}^{\top}=(\vec 0^{\top},...,\vec 0^{\top},\vec a^{\top}(0), \vec a^{\top}(1),...)$ a vector that has first $\sum_{i=1}^s K_i=K_1+...+K_s$ zero vectors $\vec 0^{\top}=(\underbrace{0,...,0}_{T})$, next vectors $\vec a(0), \vec a(1),...$ are constructed from coefficients of the functional $A\zeta$ by formula (6).

Rewrite the relation (13) in the matrix form
$$
\textbf{R} \vec {\mathbf a}=\textbf{B} \vec {\mathbf c},
$$
where  $\vec {\mathbf c}^{\top}=(\vec c^{\top}(k))_{k\in \Gamma}$ is a vector of the unknown coefficients. The linear operator $\textbf{B}$ is defined by the matrix
 $$
\textbf{B}=\begin{pmatrix}
B_{s,s} & B_{s,s-1} &\dots & B_{s,1} & B_{s,n}\\
B_{s-1,s} & B_{s-1,s-1} &\dots &B_{s-1,1} & B_{s-1,n}\\
\dots&\dots&\dots&\dots\\
B_{1,s} & B_{1,s-1} &\dots &B_{1,1} & B_{1,n}\\
B_{n,s} & B_{n,s-1} & \dots & B_{n,1} & B_{n,n}\\
\end{pmatrix},
$$
constructed with the help of the block-matrices
\begin{multline*}
B_{lm}=\left\{B_{lm}(k,j)\right\}_{k=-M_{l-1-N_l-1}}^{-M_{l}} {} _{j=-M_{m-1}-N_m-1}^{-M_{m}},\\
 B_{lm}(k,j)=B(k-j),\, l,m=1,...,s,
\end{multline*}

$$B_{ln}(k,j)=\left\{B_{ln}(k,j)\right\}_{k=-M_{l-1-N_l-1}}^{-M_{l}} {} _{j=0}^{\infty},\, B_{ln}(k,j)=B(k-j),\, l=1,...,s,$$

$$B_{nl}(k,j)=\left\{B_{nl}(k,j)\right\}_{k=0}^{\infty} {} _{j=-M_{m-1}-N_m-1}^{-M_{m}},\, B_{nl}(k,j)=B(k-j),\, m=1,...,s,$$

$$B_{nn}(k,j)=\left\{B_{nn}(k,j)\right\}_{k=0}^{\infty} {} _{j=0}^{\infty},\, B_{nn}(k,j)=B(k-j).$$

The linear operator $\textbf{R}$ is defined by the corresponding matrix, which is constructed in the same manner as matrix $\textbf{B}$.

The unknown coefficients $\vec
{c}(k),k\in \Gamma$ are determined from the equation
\begin{equation}
\label{c}
\vec {\mathbf c}=\textbf{B}^{-1}\textbf{R} \vec {\mathbf{a}},
\end{equation}
where the $k$-th component of the vector $\vec {\mathbf c}$ is the $k$-th component of vector $\textbf{B}^{-1}\textbf{R}\vec {\mathbf{a}}$:
\begin{equation}\label{ck}
\vec c(k)=(\textbf{B}^{-1}\textbf{R} \vec {\mathbf{a}})(k), \, k\in \Gamma.
\end{equation}
We will suppose that the operator $\textbf{B}$ has the inverse matrix.

The mean-square error of the optimal estimate $\widehat {A \vec \zeta}$
is calculated by the formula (9) and is of the form
$$
\Delta(\vec h, f^{\vec \zeta},f^{\vec
\theta})=E|A\vec
{\zeta}-\widehat{A\vec \zeta}|^{2}=
$$
$$
=\sum_{{\widetilde{j}=0}}^\infty \sum_{{\widetilde{k}=0}}^\infty \vec a^{\top}(\widetilde{j})\frac{1}{2\pi} \int_{-\pi}^{\pi} f^{\vec \zeta}(\lambda)(f^{\vec \zeta}(\lambda)+f^{\vec \theta}(\lambda))^{-1}  f^{\vec \theta}(\lambda)e^{-i\lambda(\widetilde{j}-\widetilde{k})} d\lambda \cdot \overline{\vec a(\widetilde{k})}+
$$
$$
+\sum_{n\in \Gamma} \sum_{k\in \Gamma} \vec c^{\top}(\widetilde{j})\frac{1}{2\pi} \int_{-\pi}^{\pi} (f^{\vec \zeta}(\lambda)+f^{\vec \theta}(\lambda))^{-1} e^{-i\lambda(n-k)} d\lambda  \cdot \overline{\vec c(k)}=
$$
\begin{equation}
\label{3.8}
=\langle{\textbf{D} \vec {\mathbf a},\vec{\mathbf a}}\rangle+\langle{\textbf{B} \vec {\mathbf c},\vec{\mathbf c}}\rangle,
\end{equation}
\noindent where
$\langle{a,b}\rangle$ denotes the scalar product, $\textbf{D}$ is defined by the corresponding matrix, which is constructed in the same manner as matrix $\textbf{B}$, with elements
$$
D(\widetilde{k}-\widetilde{j})=\frac{1}{2\pi}\int_{-\pi}^{\pi}{\left[ { f^{\vec
\zeta}(\lambda)(f^{\vec \zeta}(\lambda)+f^{\vec
\theta}(\lambda))^{-1}f^{\vec \theta}(\lambda)}
\right]}^{\top}e^{i(\widetilde{j}-\widetilde{k})\lambda} d\lambda,
$$
\[\widetilde{k}\geq0,\,\widetilde{j}\geq0.
\]

\noindent See Moklyachuk \& Masyutka (2012) for more details.

The following statement holds true.

\begin{theorem}
\label{theorem3.3}
Let $\zeta(j)$ and
$\theta(j)$ be uncorrelated T-PC stochastic sequences with the
spectral density matrices $f^{\vec \zeta}(\lambda)$ and $f^{\vec
\theta}(\lambda)$ of T-variate stationary sequences $\vec \zeta(\widetilde{j})$
and $\vec \theta(\widetilde{j})$, respectively. Assume that  $f^{\vec
\zeta}(\lambda)$ and $f^{\vec \theta}(\lambda)$ satisfy the
minimality condition(7). Assume that condition (4) is satisfied and operator $\textbf{B}$ is invertible.
The spectral characteristic $\vec h(f^{\vec \zeta},f^{\vec \theta})$ and
the mean square error $\Delta(f^{\vec \zeta},f^{\vec \theta})$ of the
optimal linear estimate of the functional $A\vec {\zeta}$
based on observations of the sequence $\vec
{\zeta}(\widetilde{j})+\vec{\theta}(\widetilde{j})$ at points  $\widetilde{j}\in \{..., -n,..., -1\}
\setminus \widetilde{S}$, are calculated by formulas (11) and (16).
\end{theorem}

Consider the mean-square estimation problem of $A \vec \zeta$ based on observations of the sequence $\vec {\zeta}(\widetilde{j})$ at points  $\widetilde{j}\in \{..., -n,..., -1\}
\setminus \widetilde{S}$. In this case the spectral density $f^{\vec \theta}(\lambda)=0$. The spectral characteristic $\vec h(f^{\vec \zeta})$ of the estimate $\widehat{A \vec \zeta}$ is of the form
\begin{equation}
\label{2.15}
\vec h^{\top}(f^{\vec \zeta})= A^{\top} (e^{i\lambda}) -
C^{\top} (e^{i\lambda }) \left[ f^{\vec
\zeta}(\lambda) \right]^{-1},
\end{equation}
where unknown coefficients $\vec c(k),\, k\in \Gamma$ are determined from the relation
  \begin{equation}\label{coef}
  \textbf{B}\vec {\mathbf c}=\vec {\mathbf a}
  \end{equation}
  or
  $$\vec {\mathbf c}=\textbf{B}^{-1}\vec {\mathbf a},$$ where
the linear operator $\textbf{B}$ is defined by the matrix
 $$
\textbf{B}=\begin{pmatrix}
B_{s,s} & B_{s,s-1} &\dots & B_{s,1} & B_{s,n}\\
B_{s-1,s} & B_{s-1,s-1} &\dots &B_{s-1,1} & B_{s-1,n}\\
\dots&\dots&\dots&\dots\\
B_{1,s} & B_{1,s-1} &\dots &B_{1,1} & B_{1,n}\\
B_{n,s} & B_{n,s-1} & \dots & B_{n,1} & B_{n,n}\\
\end{pmatrix},
$$
constructed with the help of the block-matrices
\begin{multline*}
B_{lm}=\left\{B_{lm}(k,j)\right\}_{k=-M_{l-1-N_l-1}}^{-M_{l}} {} _{j=-M_{m-1}-N_m-1}^{-M_{m}},\\
 B_{lm}(k,j)=B(k-j),\, l,m=1,...,s,
\end{multline*}
$$B_{ln}(k,j)=\left\{B_{ln}(k,j)\right\}_{k=-M_{l-1-N_l-1}}^{-M_{l}} {} _{j=0}^{\infty},\, B_{ln}(k,j)=B(k-j),\, l=1,...,s,$$

$$B_{nl}(k,j)=\left\{B_{nl}(k,j)\right\}_{k=0}^{\infty} {} _{j=-M_{m-1}-N_m-1}^{-M_{m}},\, B_{nl}(k,j)=B(k-j),\, m=1,...,s,$$

$$B_{nn}(k,j)=\left\{B_{nn}(k,j)\right\}_{k=0}^{\infty} {} _{j=0}^{\infty},\, B_{nn}(k,j)=B(k-j),$$

\noindent with elements
$$
B(k-j)=\frac{1}{2\pi}\int_{-\pi}^{\pi}{\left[
{(f^{\vec \zeta}(\lambda))^{-1}}
\right]}^{\top}e^{i(j-k)\lambda} d\lambda,
$$
\[k\in\Gamma,\,j\in \Gamma.
\]
The
mean square error $\Delta(f^{\vec \zeta})$ is defined by the formula
\begin{equation}
\label{2.16}
\Delta(f^{\vec
\zeta})=\langle{\vec {\mathbf c},\vec {\mathbf a}}\rangle.
\end{equation}

Thus, in the case without noise we have the following result.

\begin{corollary}
\label{cor3.5}
Let $\zeta(j)$ be a
T-PC stochastic sequence with the spectral density matrix $f^{\vec
\zeta}(\lambda)$ of T-variate stationary sequence $\vec \zeta(j)$.
Assume that $f^{\vec \zeta}(\lambda)$ satisfies the minimality
condition
\begin{equation}
\label{3.9}
\int_{-\pi}^{\pi}{Tr{\left[ {(f^{\vec \zeta}(\lambda))^{-1}} \right]}}
d\lambda <+{\infty}.
\end{equation}
Assume that condition (4) is satisfied and operator $\textbf{B}$ is invertible.
Then the optimal linear estimate of $A\vec \zeta$
based on observations of $\vec
{\zeta}(\widetilde{j})$ at points  $\widetilde{j}\in \{..., -n,..., -1\}
\setminus \widetilde{S}$, is given by the formula
$$
\widehat{{A}\vec\zeta}=\int_{-\pi}^{\pi}\vec h^{\top}(f^{\vec
\zeta})Z^{\vec\zeta}(d\lambda)=\int_{-\pi}^{\pi}\sum_{\nu=1}^{T}h_{\nu}(f^{\vec
\zeta})Z_{\nu}^{\vec\zeta}(d\lambda).
$$
The spectral characteristic $\vec h(f^{\vec \zeta})$ and the
mean square error $\Delta(f^{\vec \zeta})$ of $ \widehat{A\vec\zeta}$ are
calculated by formulas (17) and (19).
\end{corollary}

Let us consider the mean-square estimation problem of functional
$$A_N \zeta=\sum_{j=1}^{N\cdot T}a(j)\zeta(j)$$ that depends on unknown values of T-PC stochastic sequence $\zeta(j)$, based on observations of the sequence $ {\zeta}(j)+\theta(j)$ at points  $j\in \{..., -n,..., -1,0\}
\setminus {S}$. $\theta(j)$ is  uncorrelated with $\zeta(j)$ T-PC stochastic sequence.

Using Proposition 2.3,  the
 linear functional $A_N{\zeta}$ can be written as follows
 $$
A_N\zeta=\sum_{j=1}^{N\cdot T}{a}(j){\zeta}(j)=\sum_{{\widetilde{j}}=0}^{N-1}{a}(\widetilde{j})\sum_{\nu=1}^{T}e^{2\pi
i\widetilde{j}\nu/T}\zeta(\nu+\widetilde{j}T)=
$$
$$
=\sum_{\widetilde{j}=0}^{N-1}\sum_{\nu=1}^{T}{a}(\widetilde{j})e^{2\pi
i\widetilde{j}\nu/T}\zeta_{\nu}(j)=\sum_{\widetilde{j}=0}^{N-1}\vec {a}^{\top}(\widetilde{j})\vec
{\zeta}(\widetilde{j})=A_N\vec {\zeta},
$$

\noindent where $\vec {a}^{\top}(\widetilde{j})$ is defined by relation (6),
$\vec {\zeta}(\widetilde{j})=\left\{ \zeta_{\nu}(\widetilde{j})
\right\}_{\nu = 1}^{T}$ is $T$-variate stationary sequence, obtained by the $T$-blocking (3) of univariate $T$-PC sequence $\zeta(j), \, j \geq1$.

Let $\vec {\zeta}(j)$ and $\vec {\theta}(j)$ be uncorrelated T-variate
stationary stochastic sequences with the spectral density matrices
$f^{\vec \zeta}(\lambda)=\left\{ f^{\vec \zeta}_{\nu\mu}(\lambda)
\right\}_{\nu,\mu = 1 }^{T}$ and $f^{\vec \theta}(\lambda)=\left\{
f^{\vec \theta}_{\nu\mu}(\lambda) \right\}_{\nu,\mu= 1}^{T}$,
respectively. Consider the problem of optimal linear estimation of
the functional
\begin{equation}
\label{an}
A_N\vec {\zeta}=\sum_{\widetilde{j}=0}^{N-1}\vec {a}^{\top}(\widetilde{j})\vec
{\zeta}(\widetilde{j}),
 \end{equation}
that depends on the unknown values of sequence  $\vec {\zeta}(\widetilde{j})$,
based on observations of the sequence $\vec
{\zeta}(\widetilde{j})+\vec{\theta}(\widetilde{j})$ at points $\widetilde{j}\in\{...,-n,...,-1\}
\setminus \widetilde{S},$ $\widetilde{S}=\bigcup _{l=1}^{s}\{-M_l,\dots,-M_{l-1}-N_{l}-1\},\,
 M_l=\sum_{k=0}^l(N_k+K_k),\,\,N_0=K_0=0$.

The estimate
\begin{equation}
\label{estim}
\widehat{{A}_N\vec\zeta}=\int_{-\pi}^{\pi}\vec h_N^{\top}(e^{i\lambda})Z^{\vec\zeta}(d\lambda)
\end{equation}
of the functional $A_N \vec \zeta$ is defined by the spectral characteristic $\vec h_N(e^{i\lambda})\in L_2^s(f^{\vec \zeta}+f^{\vec \theta})$.

Denote by $\vec {\mathbf a_N}^{\top}=(\vec 0^{\top},...,\vec 0^{\top},\vec a^{\top}(0),..., \vec a^{\top}(N-1), \vec 0^{\top}, \vec 0^{\top}, ...)$ a vector that has first $ \sum_{i=1}^s K_i$ zero vectors $\vec 0^{\top}$, next $N$ vectors $\vec a(0),..., \vec a(N-1)$ are constructed from coefficients of the functional $A_N\zeta$ by formula (6).

With the help of Hilbert space projection method we can derive the following relations for all $m\in \Gamma$:
$$\sum_{\widetilde{j}=0}^{N-1} \vec a^{\top}(\widetilde{j})\frac{1}{2\pi}\int_{-\pi}^{\pi}f^{\vec \zeta}(\lambda)(f^{\vec \zeta}(\lambda)+f^{\vec \theta}(\lambda))^{-1}e^{i\lambda(\widetilde{j}-m)}d\lambda=
$$
\begin{equation}
\label{3}
\sum_{n\in \Gamma} \vec c^{\top}(n)\frac{1}{2\pi}\int_{-\pi}^{\pi}(f^{\vec \zeta}(\lambda)+f^{\vec \theta}(\lambda))^{-1}e^{i\lambda(n-m)}d\lambda.
\end{equation}

Denote by $\textbf{R}_N$ the linear operator which is defined as follows: $\textbf{R}_N(k,j)=\textbf{R}(k,j),\, j\leq N-1,$  $\textbf{R}_N(k,j)=0,\, j>N-1.$ Then we can rewrite the relations (3) in the matrix form
$$
\textbf{R}_N\vec{\mathbf a_N}=\textbf{B} \vec {\mathbf c}.
$$

The unknown vectors $\vec c(k), \, k\in \Gamma,$ are determined from the equation
$$
\vec {\mathbf c}=\textbf{B}^{-1}\textbf{R}_N\vec{\mathbf a_N}.
$$

The spectral characteristic of the optimal estimate $\widehat{{A}_N\vec\zeta}$ is calculated by formula
\begin{equation}
\label{hN}
\vec h^{\top}_N(e^{i\lambda})= \left( A^{\top}_N
(e^{i\lambda} )f^{\vec \zeta}(\lambda) - C^{\top} (e^{i\lambda
})\right)\, \left[ f^{\vec \zeta}(\lambda)+f^{\vec \theta}(\lambda)
\right]^{-1},
\end{equation}
\noindent where
$$
A_N(e^{i\lambda})=\sum_{\widetilde{j}=0}^{N-1} \vec
{a}(\widetilde{j})e^{i\widetilde{j}\lambda}.$$

The mean-square error of the optimal estimate $\widehat{{A}_N\vec\zeta}$ is calculated by formula
$$
\Delta(\vec h_N, f^{\vec \zeta},f^{\vec
\theta})=E|A_N\vec
{\zeta}-\widehat{A_N\vec \zeta}|^{2}=
$$
$$
=\sum_{{\widetilde{j}=0}}^{N-1} \sum_{{\widetilde{k}=0}}^{N-1} \vec a^{\top}(\widetilde{j})\frac{1}{2\pi} \int_{-\pi}^{\pi} f^{\vec \zeta}(\lambda)(f^{\vec \zeta}(\lambda)+f^{\vec \theta}(\lambda))^{-1}  f^{\vec \theta}(\lambda)e^{-i\lambda(\widetilde{j}-\widetilde{k})} d\lambda \cdot \overline{\vec a(\widetilde{k})}+
$$
$$
+\sum_{n\in \Gamma} \sum_{k\in \Gamma} \vec c^{\top}(\widetilde{j})\frac{1}{2\pi} \int_{-\pi}^{\pi} (f^{\vec \zeta}(\lambda)+f^{\vec \theta}(\lambda))^{-1} e^{-i\lambda(n-k)} d\lambda  \cdot \overline{\vec c(k)}=
$$
\begin{equation}
\label{deltaN}
=\langle{\textbf{D}_N \vec {\mathbf a}_N,\vec{\mathbf a}_N}\rangle+\langle{\textbf{B} \vec {\mathbf c},\vec{\mathbf c}}\rangle,
\end{equation}
\noindent where linear operator $\textbf{D}$ is defined  as follows: $\textbf{D}_N(k,j)=\textbf{D}(k,j),\, k, j\leq N-1,$  $\textbf{D}_N(k,j)=0$ if $ k>N-1$ or $ j>N-1.$

\begin{theorem}
\label{theoremAn}
Let $\zeta(j)$ and
$\theta(j)$ be uncorrelated T-PC stochastic sequences with the
spectral density matrices $f^{\vec \zeta}(\lambda)$ and $f^{\vec
\theta}(\lambda)$ of T-variate stationary sequences $\vec \zeta(\widetilde{j})$
and $\vec \theta(\widetilde{j})$, respectively. Assume that  $f^{\vec
\zeta}(\lambda)$ and $f^{\vec \theta}(\lambda)$ satisfy the
minimality condition (7). Assume that operator $\textbf{B}$ is invertible.
The spectral characteristic $\vec h_N(e^{i\lambda})$ and
the mean square error $\Delta(\vec h_N;f^{\vec \zeta},f^{\vec \theta})$ of the
optimal linear estimate of the functional $A_N\vec {\zeta}$
based on observations of the sequence $\vec
{\zeta}(\widetilde{j})+\vec{\theta}(\widetilde{j})$ at points  $\widetilde{j}\in \{..., -n,..., -1\}
\setminus \widetilde{S}$, are calculated by formulas (24) and (25).
\end{theorem}

In the case of observation without noise we have the following result.

\begin{corollary}
\label{cor2}
Let $\zeta(j)$ be a
T-PC stochastic sequence with the spectral density matrix $f^{\vec
\zeta}(\lambda)$ of T-variate stationary sequence $\vec \zeta(j)$.
Assume that $f^{\vec \zeta}(\lambda)$ satisfies the minimality
condition (20). Assume that operator $\textbf{B}$ is invertible.
The spectral characteristic $\vec h_N(e^{i\lambda})$ and the
mean square error $\Delta(f^{\vec \zeta})$ of $ \widehat{A_N\vec\zeta}$ are
calculated by formulas
\begin{equation}
\label{hnn}
\vec h^{\top}_N(e^{i\lambda})= A^{\top}_N (e^{i\lambda}) -
C^{\top} (e^{i\lambda }) \left[ f^{\vec
\zeta}(\lambda) \right]^{-1},
\end{equation}
\begin{equation}
\label{deltann}
\Delta(f^{\vec
\zeta})=\langle{\vec {\mathbf c},\vec {\mathbf a}_N}\rangle.
\end{equation}
The linear operator $\textbf{B}$ is defined in Corollary 1, vector $\vec{\mathbf c}$ is defined by the equation $\vec{\mathbf c}= \textbf{B}^{-1}\vec {\mathbf a}_N$.
\end{corollary}

\begin{exm}
Let  $\zeta(n),\, n\in \mathbb{Z}$, be a 2-PC stochastic sequence such that $\zeta(2n)=\eta(n)$ is a univariate white noise with the spectral density $f(\lambda)=1$ and $\zeta(2n+1)=\gamma(n)$ is an uncorrelated with $\eta(n)$ univariate stationary Ornstein-Uhlenbeck sequence with the spectral density
$g(\lambda)=\frac{1}{|1-e^{i\lambda}|^{2}}$.
Consider the problem of estimation of the functional
$$
A_{1}\zeta=\zeta(1)+\zeta(2)
$$
based on observations of $\zeta(n),n\in \{...,-1,0\}\setminus\{-3,-2\}=\{...,-5,-4,-1,0\}$. Here $S=\{-3,-2\},\, N_1=K_1=1,\, M_1=2$.

Rewrite functional $A_1\zeta$ in the form (21)
$$
A_{1}\zeta=\zeta(1)+\zeta(2)=(1,1)\cdot
\begin{pmatrix}
\zeta_1(0)\\
\zeta_2(0)
\end{pmatrix}=\vec a^{\top}(0) \vec \zeta(0)=A_1\vec \zeta,
$$
where $\vec a(0)=(a(1+0\cdot2)e^{2\pi i1\cdot0/2},a(2+0\cdot2)e^{2\pi i2\cdot0/2})^{\top}=(1,1)^{\top},\,$ $\vec \zeta(0)=(\zeta(1+0\cdot2),\zeta(2+0\cdot2))^{\top}=(\zeta_1(0),\zeta_2(0))^{\top},$ $\tilde{S}=\{-2\}.$
The spectral density matrix of 2-variate stationary sequence  $\vec{\zeta}(n)$ is of the form
$$
f^{\vec
\zeta}(\lambda)=\begin{pmatrix} f(\lambda) & 0\\
0 & g(\lambda)
\end{pmatrix}
$$
The matrix ${ [{f^{\vec \zeta}(\lambda)}]}^{-1}$ is of the form
$$
[{f^{\vec \zeta}(\lambda)}]^{-1}=\begin{pmatrix}
1&0\\
0&2
\end{pmatrix}+\begin{pmatrix}
0 & 0\\
0 & -1
\end{pmatrix} e^{-i\lambda}+\begin{pmatrix}
0 & 0\\
0 & -1
\end{pmatrix} e^{i\lambda}=B(0)+B(-1)e^{-i\lambda}+B(1)e^{i\lambda} $$
and satisfies the minimality
condition (20). In the last equality matrices
$$
B(0)=\begin{pmatrix}
1 & 0\\
0 & 2
\end{pmatrix} ,\, B(-1)=B(1)=\begin{pmatrix}
0 & 0\\
0 & -1
\end{pmatrix}
$$
are the Fourier coefficients of the function ${ [{f^{\vec \zeta}(\lambda)}]}^{-1}$. In order to find the spectral characteristic $\vec h_1(e^{i\lambda})$ and the mean-square error $\Delta(f^{\vec \zeta})$ of the estimate $\widehat {A_1 \vec \zeta}$ let us use the Corollary 2. To find the unknown coefficients
$$
\vec c (k)=(\textbf{B}^{-1}\vec {\mathbf a}_N)(k),\, k\in \Gamma=\widetilde{S}\cup\{0,1,...\}=\{-2,0,1,...\}
$$
we use the equation (18), where vectors $\vec {\mathbf c}^{\top}=(\vec c^{\top}(-2),\vec c^{\top}(0), \vec c^{\top}(1),...),$  $\vec {\mathbf a}^{\top}_1=(\vec 0^{\top},\vec a^{\top}(0), \vec 0^{\top},...).$ The operator $\textbf{B}$ is defined by matrix
$$
\textbf{B}=\begin{pmatrix}
B_{11} & B_{1n}\\
B_{n1} & B_{nn}
\end{pmatrix},$$
with block-matrices
$$B_{11}=\left\{B_{11}(k,j)\right\}_{k=-2} {}_{j=-2}=B(0),$$
$$
B_{1n}=\left\{B_{1n}(k,j)\right\}_{k=-2} {}_{j=0}^\infty=(B(-2)\, B(-3)\, B(-4)\,...)=(O_2\, O_2\, O_2\,...),
$$
$$
B_{n1}=\left\{B_{n1}(k,j)\right\}_{k=0} {}_{j=-2}=(B(2)\, B(3)\, B(4)\,...)^{\top}=(O_2\, O_2\, O_2\,...)^{\top},
$$
$$
B_{nn}=\left\{B_{nn}(k,j)\right\}_{k=0}^\infty {}_{j=0}^\infty=\begin{pmatrix}
B(0) & B(-1) & O_2 &...\\
B(1) & B(0) & B(-1) &...\\
O_2 & B(1) & B(0) & ...\\
. & .& .& ...
\end{pmatrix},$$
where $O_2=\begin{pmatrix} 0&0\\
0&0
\end{pmatrix}$.

The inverse matrix $\textbf{B}^{-1}$ can be represented in the form
$$
\textbf{B}^{-1}=\begin{pmatrix}
B_{11}^{-1} & 0\\
0 & B_{nn}^{-1}
\end{pmatrix},
$$
where $B_{11}^{-1}=(B(0))^{-1}$, $B_{nn}^{-1}$ is the inverse matrix to $B_{nn} $. To find $B_{nn}^{-1}$ we use that matrix $[f^{\vec \zeta}(\lambda)]^{-1}$ admits factorization
$$
[f^{\vec \zeta}(\lambda)]^{-1}=\sum_{j=-\infty}^{\infty}B(j)e^{ij\lambda}=\left(\sum_{k=0}^\infty \psi(k)e^{-ik\lambda}\right)\left(\sum_{k=0}^\infty \psi(k)e^{-ik\lambda}\right)^*=
$$
$$
=\left(\left(\sum_{k=0}^\infty \varphi(k)e^{-ik\lambda}\right)^* \left(\sum_{k=0}^\infty \varphi(k)e^{-ik\lambda}\right)\right)^{-1},
$$
where $\psi(0)=\begin{pmatrix}
1&0\\
0&1
\end{pmatrix},$  $\psi(1)=\begin{pmatrix}
0&0\\
0&-1
\end{pmatrix},$ $\psi(k)=O_2, k\geq2$
   and  $\varphi(0)=\begin{pmatrix}
1&0\\
0&1
\end{pmatrix},$ $\varphi(k)=\begin{pmatrix}
0&0\\
0&1
\end{pmatrix}, k\geq1$.

If we denote by $\Psi$ and $\Phi$ linear operators determined by matrices with elements $\Psi(i,j)=\psi(j-i),$ $\Phi(i,j)=\varphi(j-i),$ for $0\leq i\leq j,$ $\Psi(i,j)=0,$ $\Phi(i,j)=0,$ for $0\leq j< i.$ Then elements of the matrix $B_{nn}$ can be represented in the form $B_{nn}(i,j)=(\Psi\Psi^*)(i,j)$. It is not hard to verify that $\Psi\Phi=\Phi\Psi=I.$ This makes possible to write elements of $B_{nn}^{-1}$ in the form $B_{nn}^{-1}(i,j)=(\Phi^*\Phi)(i,j)=\sum_{l=0}^{\min(i,j)}(\varphi(i-l))^*\varphi(j-l).$

Using equation $\vec{\mathbf c}= \textbf{B}^{-1}\vec {\mathbf a}_N$ we can represent the unknown coefficients $\vec c(k), k\in \Gamma$ in the form
$$
\vec c(-2)=\vec 0,
$$
$$
\vec c(0)=B_{nn}^{-1}(0,0)\vec a(0),
$$
$$
\vec c(1)=B_{nn}^{-1}(1,0)\vec a(0),$$
$$...$$
$$
\vec c(i)=B_{nn}^{-1}(i,0)\vec a(0),\, i\geq2.
$$

The spectral characteristic $\vec h_1(e^{i\lambda})$ is determined by the formula (26)
$$
\vec h_1^{\top}(e^{i\lambda})=-\vec c^{\top}(0)B(-1)e^{-i\lambda}=-B_{nn}^{-1}(0,0) \vec a(0) B(-1)e^{-i\lambda}.
$$
Since $B_{nn}^{-1}=\varphi^*(0)\varphi(0)=\begin{pmatrix}
1&0\\
0&1
\end{pmatrix}$, the spectral characteristic is of the form
$$
\vec h_1^{\top}(e^{i\lambda})=-(0,-1)e^{-i\lambda}.
$$

The optimal linear estimate $\widehat{A_1\vec \zeta}$ can be calculated by the formula (22)
$$\widehat{ A_1\vec \zeta}=\zeta_2(-1)=\zeta(0).$$

The mean-square error of the estimate $\widehat{A_1\vec \zeta}$ determined by (27) equals
$$\Delta(f^{\vec \zeta})=\langle\vec {\mathbf c},\vec {\mathbf a}_1\rangle=2. $$
\end{exm}

\section{Minimax (robust) method of linear extrapolation problem}

Let $f(\lambda)$ and $g(\lambda)$ be the spectral density
matrices of  $T$-variate stationary
sequences $\vec \zeta(j)$ and $\vec \theta(j)$, obtained by $T$-blocking (3) of $T$-PC sequences  $\zeta(j)$ and $\theta(j)$, respectively.

The obtained formulas may be applied for finding the spectral characteristic and
the mean square error of the optimal linear estimate of the functionals $A \vec \zeta$ and $A_N \vec \zeta$
only under the condition that the spectral density
matrices $f(\lambda)$ and $g(\lambda)$ are exactly known.
If the density matrices are not known exactly while a set
$D=D_{f}\times D_{g}$ of possible spectral densities is given, the
minimax (robust) approach to estimation of functionals from unknown
values of stationary sequences is reasonable. In this case we find the
estimate which minimizes the mean square error for all spectral
densities from the given set simultaneously.

\begin{definition} For a given class of
pairs of spectral densities $D=D_{f}\times D_{g}$ the spectral density
matrices $f^{0}(\lambda)\in D_{f}$, $g^{0}(\lambda)\in D_{g}$ are
called  {the least favorable} in $D$ for the optimal linear
estimation of the functional $A\vec\zeta$ if
$$
\Delta(f^0,g^0)=\Delta(\vec h(f^0,g^0);f^0,g^0)={\max_{\substack
{(f,g)\in D}}} \Delta (\vec h(f,g);f,g).
$$
\end{definition}

\begin{definition}
For a given class of
pairs of spectral densities $D=D_{f}\times D_{g}$ the spectral
characteristic $\vec h^0(\lambda)$ of the optimal linear estimate of
the functional $A\vec \zeta$ is called minimax (robust) if
$$
\vec h^0(\lambda)\in H_{D}=\mathop {\bigcap }\limits_{(f,g)\in D}
L_{2}^{s}(f+g),$$
$$\mathop {\min_{\substack {\vec h\in H_{D}}}} {\max_
{\substack {(f,g)\in D}}} \Delta(\vec h;f,g)=\mathop {\max_{\substack
{(f,g)\in D}}} \Delta(\vec h^0;f,g).
$$
\end{definition}

Taking into consideration these definitions and the obtained relations we can
verify that the following lemmas hold true.

\begin{lemma}
\label{lem4.1}
The spectral density matrices
$f^{0}(\lambda)\in D_{f}$, $g^{0}(\lambda)\in D_{g}$, that satisfy the minimality
 condition (7), are the least favorable in the class D for the optimal
linear  estimation of $A\vec\zeta$, if the Fourier coefficients of
the matrix functions
$$
(f^0(\lambda)+g^0(\lambda))^{-1},\quad
f^0(\lambda)(f^0(\lambda)+g^0(\lambda))^{-1},
$$
$$
f^0(\lambda)(f^0(\lambda)+g^0(\lambda))^{-1}g^0(\lambda)
$$
define matrices $\textbf{B}^0, \textbf{R}^0, \textbf{D}^0$, that determine a solution of the constrained optimization problem
$$
{\max_{\substack {(f,g)\in
D}}}(\langle{\textbf{R}\vec{a}},\textbf{B}^{-1}\textbf{R}\vec{a}\rangle)+\langle{\textbf{D}\vec{a},\vec{a}}\rangle)=\langle{\textbf{R}^0\vec{a}},(\textbf{B}^0)^{-1}\textbf{R}^0\vec{a}\rangle)+\langle{\textbf{D}^0\vec{a},\vec{a}}\rangle.
$$
The minimax spectral characteristic
$\vec h^0=\vec h(f^0,g^0)$ is given by (11), if $\vec h(f^0,g^0)\in H_{D}$.
\end{lemma}

In the case of observations of the sequence without noise the following corollary holds true.

\begin{corollary}
\label{lem4.2}
The spectral density matrix
$f^{0}(\lambda)\in D_{f}$, that satisfies the minimality  condition  (20), is
the least favorable in the class $D_{f}$ for the optimal linear estimation
of $A\vec\zeta$
based on observations of $\vec
{\zeta}(\widetilde{j})$ at points  $\widetilde{j}\in \{..., -n,..., -1\}
\setminus \widetilde{S}$, if the Fourier
coefficients of the matrix function $(f^0(\lambda))^{-1}$ define the matrix
$\textbf{B}^0$, that determine a solution of the constrained optimization problem
$$
\max_{\substack {f\in D_f}}
\langle{\textbf{B}^{-1}\vec{a},\vec{a}}\rangle=\langle{(\textbf{B}^0)^{-1}\vec{a},\vec{a}}\rangle.
$$
The minimax spectral characteristic $\vec h^0=\vec h(f^0)$
is given by (17), if $\vec h(f^0)\in H_{D}$.
\end{corollary}

The least favorable spectral densities $f^{0}(\lambda)\in D_{f}$,
$g^{0}(\lambda)\in D_{g}$ and the minimax spectral characteristic
$\vec h^0=\vec  h(f^0,g^0)$ form a saddle point of the function $\Delta(\vec h;f,g)$
on the set $H_{D}\times D$. The saddle point inequalities
$$
\Delta(\vec h^0;f,g)\leq\Delta(\vec h^0;f^0,g^0)\leq\Delta(\vec h;f^0,g^0), \quad \forall
\vec h\in H_{D}, \forall f\in D_{f}, \forall g\in D_{g}
$$
\noindent hold true when $\vec h^0=\vec h(f^0,g^0)$, $\vec h(f^0,g^0)\in H_{D}$ and
$(f^0,g^0)$ is a solution of the constrained optimization problem
\begin{equation}
 \label{copt1}
\Delta\left(\vec h(f^0,g^0);f,g\right)\rightarrow \sup, \,\, (f,g)\in D_f\times D_g.
\end{equation}

\noindent
The linear functional $\Delta(\vec h(f^0,g^0);f,g)$
is calculated by the formula

\begin{multline*}
\nonumber
\Delta(\vec h(f^0,g^0);f,g)=\frac{1}{2\pi}\int_{-\pi}^{\pi}\left({A}^{\top}(e^{i\lambda})g^0(\lambda)+
({C}^0(e^{i\lambda}))^{\top}\right)
\times
\\
\times
(f^0(\lambda)+g^0(\lambda))^{-1}
f(\lambda)(f^0(\lambda)+g^0(\lambda))^{-1}
\bigg({A}^{\top}(e^{i\lambda})g^0(\lambda)+
({C}^0(e^{i\lambda}))^{\top}\bigg)^*d\lambda+\\
+
 \frac{1}{2\pi}\int_{-\pi}^{\pi}\left({A}^{\top}(e^{i\lambda})f^0(\lambda)-
({C}^0(e^{i\lambda}))^{\top}\right)(f^0(\lambda)+g^0(\lambda))^{-1}\times
\\
\times g(\lambda)(f^0(\lambda)+g^0(\lambda))^{-1}
\left({ A}^{\top}(e^{i\lambda})f^0(\lambda)-
({C}^0(e^{i\lambda}))^{\top}\right)^*d\lambda,
\end{multline*}

\noindent  where ${C}^0(e^{i\lambda})=\sum_{n\in \Gamma} \vec c^{~0}(n) e^{in\lambda},$ column vectors $\vec c^{~0}(n)=((\textbf{B}^0)^{-1} \textbf{R}^0 \vec{\mathbf a})(n)$.

The constrained optimization problem (28) is equivalent to the unconstrained optimization problem (see Pshenichnyi, 1971):
\begin{equation}
\label{copt2}
\Delta_D(f,g)=-\Delta(\vec h(f^0,g^0);f,g)+\delta((f,g)\left|D_f\times D_g\right.)\rightarrow \inf,
\end{equation}
where $\delta((f,g)|D_f\times D_g)$ is the indicator function of the set $D=D_f\times D_g$.
A solution of the problem (29) is characterized by the condition $0 \in \partial\Delta_D(f^0,g^0),$ where $\partial\Delta_D(f^0,g^0)$ is the subdifferential of the convex functional $\Delta_D(f,g)$ at point $(f^0,g^0)$, see Rockafellar (1997).

The form of the functional $\Delta(\vec h(f^0,g^0);f,g)$  admits finding the derivatives and differentials of the functional in the space $L_1\times L_1$. Therefore the complexity  of the optimization problem (29) is determined by the complexity of calculating of subdifferentials of the indicator functions  $\delta((f,g)|D_f\times D_g)$  of the sets $D_f\times D_g$
(see Ioffe \&  Tihomirov, 1979)

Taking into consideration the introduced definitions and the derived relations we can verify that the following lemma holds true.

\begin{lemma}
\label{lem4.3}
Let $(f^0,g^0)$ be a solution to the optimization problem (29). The spectral densities  $f^0(\lambda)$, $g^0(\lambda)$
are the least favorable in the class $D=D_f\times D_g$ and the spectral characteristic  $\vec h^0= \vec h(f^0,g^0)$
is the minimax of the optimal linear estimate of the functional  $A\vec{\zeta}$  if  $\vec h(f^0,g^0) \in H_D$.
\end{lemma}

In the case of estimation of the functional based on observations without noise we have the following statement.
\begin{lemma}
\label{lem4.4}
Let $f^0(\lambda)$ satisfies the
condition (20) and be a solution of the constrained optimization problem
\begin{equation}
\label{4.17}
\Delta(\vec h(f^0);f)\rightarrow{sup},
f(\lambda)\in D_{f},
\end{equation}
\[
\Delta(\vec h(f^0);f)=\frac{1}{2\pi}\int_{-\pi}^{\pi}\left({C}^0(e^{i\lambda})\right)^{\top}(f^0(\lambda))^{-1}f(\lambda)(f^0(\lambda))^{-1}
\overline{\left({C}^0(e^{i\lambda})\right)}d\lambda,
\]
 where ${C}^0(e^{i\lambda})=\sum_{n\in \Gamma} \vec c^{~0}(n) e^{in\lambda},$  column vectors $\vec c^{~0}(n)=((\textbf{B}^0)^{-1}  \vec {\mathbf a})(n)$.

Then $f^0(\lambda)$ is the least favorable
spectral density matrix for the optimal linear estimation of
$A\vec\zeta$
based on observations of $\vec
{\zeta}(\widetilde{j})$ at points  $\widetilde{j}\in \{..., -n,..., -1\}
\setminus \widetilde{S}$.
The minimax spectral
characteristic $\vec h^0=\vec h(f^0)$ is given by (17), if $\vec h(f^0)\in
H_{D}$.
\end{lemma}

\section{The least favorable spectral densities in the class $D=D_{0}\times D_V^U$ }

Let $f(\lambda)$ and $g(\lambda)$ be the spectral density
matrices of  $T$-variate stationary
sequences $\vec \zeta(j)$ and $\vec \theta(j)$, obtained by $T$-blocking (3) of $T$-PC sequences  $\zeta(j)$ and $\theta(j)$, respectively.

Consider the problem of minimax estimation of the functional
$A\vec\zeta$ based on observations of the sequence $\vec
{\zeta}(\widetilde{j})+\vec{\theta}(\widetilde{j})$ at points  $\widetilde{j}\in \{..., -n,..., -1\}
\setminus \widetilde{S}$, under the condition that the spectral
density matrices $f(\lambda)$ and $g(\lambda)$  belong to the class $D=D_{0}\times D_V^U$, where
$$
D_0^{1}={\left\{ {f(\lambda)|\,\frac{1}{2\pi}\int_{-\pi}^{\pi}f(\lambda)d\lambda=P}
\right\}},
$$
$$
D_V^{U1}={\left\{ {g(\lambda)|\,V(\lambda)\leq g(\lambda)\leq U(\lambda), \frac{1}{2\pi}\int_{-\pi}^{\pi}g(\lambda)d\lambda=Q}
\right\}},
$$

$$
D_0^{2}={\left\{ {f(\lambda)|\,\frac{1}{2\pi}\int_{-\pi}^{\pi}\textmd{Tr } f(\lambda)d\lambda=p}
\right\}},
$$
$$
D_V^{U2}={\left\{ {g(\lambda)|\,\textmd{Tr } V(\lambda)\leq \textmd{Tr } g(\lambda)\leq \textmd{Tr } U(\lambda), \frac{1}{2\pi}\int_{-\pi}^{\pi} \textmd{Tr } g(\lambda)d\lambda=q}
\right\}},
$$

\noindent where  $P, Q$ are known positive definite Hermitian matrices, spectral densities $V(\lambda), U(\lambda)$ are known and fixed, $p, q$ are known and fixed numbers.

With the help of the method of Lagrange multipliers we can find that solution
$(f^{0}(\lambda), g^0(\lambda))$ of the constrained optimization problem (28) satisfy
the following relations for these sets of admissible spectral densities.

For the pair $D_{0}^1\times D_V^{U1}$ we have relations
\begin{equation}
\label{1g}
(g^0(\lambda)\overline{ A(e^{i\lambda})}+\overline{C^0(e^{i\lambda})})((g^0(\lambda))^{\top} A(e^{i\lambda})+C^0(e^{i\lambda}))^{\top}=
\end{equation}
\[(f^0(\lambda)+g^0(\lambda))\overline{\vec \alpha}\vec \alpha^{\top}(f^0(\lambda)+g^0(\lambda)),\]
\begin{equation}
\label{1f}
(f^0(\lambda) \overline{A(e^{i\lambda})}-\overline{C^0(e^{i\lambda})})((f^0(\lambda))^{\top}  A(e^{i\lambda})-C^0(e^{i\lambda}))^{\top}=
\end{equation}
\[(f^0(\lambda)+g^0(\lambda))(\overline{\vec \beta}\vec \beta^{\top}+\psi_1(\lambda)+\psi_2(\lambda))(f^0(\lambda)+g^0(\lambda)),\]
\noindent where $\vec \alpha,\vec \beta$ are
Lagrange multipliers, $\psi_1(\lambda)\leq0$ and $\psi_1(\lambda)=0$ if $g^0(\lambda)\geq V(\lambda)$, $\psi_2(\lambda)\geq0$ and $\psi_2(\lambda)=0$ if $g^0(\lambda)\leq U(\lambda)$.

For the pair $D_{0}^2\times D_V^{U2}$ we have relations
\begin{equation}
\label{2g}
(g^0(\lambda)\overline{ A(e^{i\lambda})}+\overline{C^0(e^{i\lambda})})((g^0(\lambda))^{\top} A(e^{i\lambda})+C^0(e^{i\lambda}))^{\top}=
\end{equation}
\[\alpha^2(f^0(\lambda)+g^0(\lambda))^2,\]
\begin{equation}
\label{2f}
(f^0(\lambda) \overline{A(e^{i\lambda})}-\overline{C^0(e^{i\lambda})})((f^0(\lambda))^{\top}  A(e^{i\lambda})-C^0(e^{i\lambda}))^{\top}=
\end{equation}
\[(\beta^2+\varphi_1(\lambda)+\varphi_2(\lambda))(f^0(\lambda)+g^0(\lambda))^2,\]
\noindent where $\alpha^2, \beta^2$ are
Lagrange multipliers, $\varphi_1(\lambda)\leq0$ and $\varphi_1(\lambda)=0$ if $\textmd{Tr } g^0(\lambda)\geq \textmd{Tr } V(\lambda)$, $\varphi_2(\lambda)\geq0$ and $\varphi_2(\lambda)=0$ if $\textmd{Tr } g^0(\lambda)\leq \textmd{Tr } U(\lambda)$.

Hence the following theorem  holds true.

\begin{theorem}
\label{theorem3}
Let the spectral densities $f^0(\lambda)$ and $g^0(\lambda)$ satisfy the minimality condition (7). The least favorable spectral densities $f^0(\lambda)$, $g^0(\lambda)$ in the class $D_{0}^1\times D_V^{U1}$ for the optimal linear extrapolation  of the functional $A\vec\zeta$ are determined by relations (31), (32).  The least favorable spectral densities $f^0(\lambda)$, $g^0(\lambda)$ in the class $D_{0}^2\times D_V^{U2}$ for the optimal linear extrapolation  of the functional $A\vec\zeta$ are determined by relations (33), (34). The minimax spectral characteristic of the optimal estimate of the functional $A\vec\zeta$ is determined by the formula (11).
\end{theorem}

In the case of observations of the sequence without noise the following corollaries hold true.

\begin{corollary}
\label{corr1}
Let the spectral density $f^0(\lambda)$  satisfies the minimality condition (20).
The least favorable spectral density $f^0(\lambda)$ in the class $D_{0}^1$ or $D_{0}^2$ for the optimal linear extrapolation  of the functional $A\vec\zeta$
based on observations of $\vec
{\zeta}(\widetilde{j})$ at points  $\widetilde{j}\in \{..., -n,..., -1\}
\setminus \widetilde{S}$ is determined by relations, respectively
\begin{equation}
\label{1gg}
(\overline{C^0(e^{i\lambda})})(C^0(e^{i\lambda}))^{\top}=f^0(\lambda)\overline{\vec \alpha}\vec \alpha^{\top}f^0(\lambda),\end{equation}

\begin{equation}
\label{2gg}
(\overline{C^0(e^{i\lambda})})(C^0(e^{i\lambda}))^{\top}=\alpha^2(f^0(\lambda))^2,\end{equation}
by the constrained optimization problem
(30) and  restrictions on the density from the corresponding class $D_{0}^1$ or $D_{0}^2$. The minimax spectral characteristic of the optimal estimate of the functional $A\vec\zeta$ is determined by the formula (17).
\end{corollary}

\begin{corollary}
\label{corr2}
Let the spectral density $f^0(\lambda)$  satisfies the minimality condition (20).
The least favorable spectral density $f^0(\lambda)$ in the class $D_V^{U1}$ or $D_V^{U2}$ for the optimal linear extrapolation  of the functional $A\vec\zeta$
based on observations of $\vec
{\zeta}(\widetilde{j})$ at points  $\widetilde{j}\in \{..., -n,..., -1\}
\setminus \widetilde{S}$ is determined by relations, respectively
\begin{equation}
\label{1ff}
(\overline{C^0(e^{i\lambda})})(C^0(e^{i\lambda}))^{\top}=f^0(\lambda)(\overline{\vec \beta}\vec \beta^{\top}+\psi_1(\lambda)+\psi_2(\lambda))f^0(\lambda),
\end{equation}
\begin{equation}
\label{2ff}
(\overline{C^0(e^{i\lambda})})(C^0(e^{i\lambda}))^{\top}=(\beta^2+\varphi_1(\lambda)+\varphi_2(\lambda))(f^0(\lambda))^2,\end{equation}
by the constrained optimization problem
(30) and  restrictions on the density from the corresponding class $D_V^{U1}$ or $D_V^{U2}$. The minimax spectral characteristic of the optimal estimate of the functional $A\vec\zeta$ is determined by the formula (17).
\end{corollary}

\section{Conclusions}

In this article we study the extrapolation of the functionals $A\zeta$ and  $A_N\zeta$  which depend on the
unobserved values of a periodically correlated stochastic sequence
${\zeta}(j)$.
Estimates are based on observations of a periodically correlated stochastic sequence
${\zeta}(j)+{\theta}(j)$ with  missing observations, that means that observations of ${\zeta}(j)+{\theta}(j)$ are known at points $j\in{\mathbb Z}\setminus S$, $j\in\{...,-n,...,-2,-1,0\}\setminus S$, $S=\bigcup _{l=1}^{s-1}\{-M_l\cdot T+1,\dots,-M_{l-1}\cdot T-N_{l}\cdot T\}$. The sequence ${\theta}(j)$ is an uncorrelated with ${\zeta}(j)$ additive noise.

The extrapolation problem is considered under the condition of spectral certainty and under the condition of spectral uncertainty. In the first case of spectral certainty  the spectral density matrices
$f(\lambda)$ and $g(\lambda)$ of the
$T$-variate stationary sequences $\vec \zeta(n)$ and $\vec \theta(n)$, obtained by $T$-blocking of $T$-PC sequences  $\zeta(j)$ and $\theta(j)$, respectively, are suppose to be known exactly.
With the help of Hilbert space projection method formulas for calculating the spectral characteristic and the mean-square error of the optimal estimate of the functionals are proposed.
In the second case of spectral uncertainty the spectral density matrices are
not exactly known while a class $D=D_{f} \times D_{g}$ of admissible
spectral densities is given. Using the minimax (robust) estimation method we derived relations that determine the least favorable spectral
densities and the minimax spectral characteristic of the optimal estimate of the functional
$A\zeta$. The problem is investigated in details for two special classes of admissible spectral densities.
In each of cases of spectral certainty and uncertainty   the case of observations of the sequence without noise ${\theta}(j)$ are presented.

 \end{document}